\documentclass[3p]{elsarticle}
\usepackage{graphicx}
\usepackage{amsmath,amssymb}
\usepackage{upgreek}
\usepackage{epstopdf}

\newcommand{\bx}{{\bf x}}
\newcommand{\by}{{\bf y}}
\newcommand{\bn}{{\bf n}}
\newcommand{\bmfact}{\eta}
\newcommand{\bphi}{\boldsymbol \phi}
\newcommand{\bG}{{\bf G}}
\newcommand{\bH}{{\bf H}}
\newcommand{\bHp}{{\bf H'}}
\newcommand{\bE}{{\bf E}}
\newcommand{\bv}{{\bf v}}
\newcommand{\de}{\mathrm{d}}
\newcommand{\bphinc}{{\boldsymbol \phi}_{\text{inc}}}
\newcommand{\beps}{{\boldsymbol \varepsilon}}
\newcommand{\I}{{\mathrm{i}}}
\newcommand{\order}{\mathcal{O}}
\newcommand{\inc}{\text{inc}}

\begin{document}
\begin{frontmatter}{}
  \title{About the Burton-Miller factor in the low frequency region}
  \author[ARI]{Wolfgang Kreuzer\corref{CORRAUTH}}
  \ead{wolfgang.kreuzer@oeaw.ac.at}
  
  \affiliation[ARI]{
    organization = {Austrian Academy of Sciences, Acoustics Research Institute},
    addressline = {Dominikanerbastei 16, 3. Stock},
    postcode = {1010},
    city = {Vienna},
    country = {Austria}
}

\cortext[CORRAUTH]{Corresponding author.}

\begin{abstract}
  The Burton-Miller method is a widely used approach in acoustics to enhance the stability of the boundary element method for exterior Helmholtz problems at so-called critical frequencies. This method depends on a coupling parameter $\bmfact$ and it can be shown that as long as $\bmfact$ has an imaginary part different from 0, the boundary integral formulation for the Helmholtz equation has a unique solution at all frequencies. A popular choice for this parameter is $\eta = \frac{\I}{k}$, where $k$ is the wavenumber. It can be shown that this choice is quasi optimal, at least in the high frequency limit. However, especially in the low frequency region, where the critical frequencies are still sparsely distributed, different choices for this factor result in a smaller condition number and a smaller error of the solution. In this work, alternative choices for this factor are compared based on numerical experiments. Additionally, a way to enhance the Burton-Miller solution with $\eta = \frac{\I}{k}$ for a sound hard scatterer in the low frequency region by an additional step of a modified Richardson iteration is introduced.
\end{abstract}
\journal{Engineering Analysis with Boundary Elements}
\end{frontmatter}{}
\section{Introduction}\label{Sec:Introduction}
The boundary element method (BEM, \cite{MarNol08V,SauSch10}) is an often used numerical method to solve exterior scattering problems in acoustics. It is a well known fact that although the exterior Helmholtz equation has a unique solution, there are certain critical frequencies where the boundary integral equation (BIE) used by the BEM has not \cite{BurMil71,Schenck68}. As the number of these critical frequencies below a given wavenumber $k$ is proportional to $k^3$ (cf. \cite{AmiHar90}), this non-uniqueness poses a major problem for the Helmholtz BEM.

One way to ensure a unique solution is the Burton-Miller method \cite{BurMil71}, that uses a coupling between the BIE
\begin{equation}\label{Equ:BIE}
  \frac{\phi(\bx)}{2} - \int_\Gamma H(\bx,\by) \phi(\by) \de\by  + \int_{\Gamma} G(\bx,\by) v(\by) \de \by = \phi_{\text{inc}}(\bx)
\end{equation}
and its derivative with respect to the normal vector $\bn_x$ 
\begin{equation}\label{Equ:dBIE}
  \frac{v(\bx)}{2} - \int_\Gamma E(\bx,\by) \phi(\by) \de\by  + \int_{\Gamma} H'(\bx,\by) v(\by) \de \by = v_{\inc}(\bx).
\end{equation}
Eq.~(\ref{Equ:dBIE}) will be called dBIE in the rest of the paper.  $\phi(\bx)$ denotes the velocity potential at a point $\bx$ on a smooth part of the boundary $\Gamma$, $v(\bx) := \frac{\partial \phi}{\partial \bn_x}\left(\bx\right) := \nabla \phi(\bx) \cdot \bn_x$ and 
\begin{align*}
  G(\bx,\by) &= \frac{e^{\I k||\bx - \by||}}{4\pi ||\bx - \by||},
               \phantom{i}H(\bx,\by) = \frac{\partial G}{\partial \bn_y}(\bx,\by) = \frac{\partial G}{\partial \by} (\bx,\by) \cdot \bn_y,\\
  H'(\bx,\by) &= \frac{\partial G}{\partial \bn_x}(\bx,\by), \quad E(\bx,\by) = \frac{\partial^2 G}{\partial \bn_x  \partial \bn_y}(\bx,\by),
\end{align*}
where $\bn_y$ is the normal vector to the boundary $\Gamma$ at a point $\by \in \Gamma$. In Eqs.~(\ref{Equ:BIE}) and (\ref{Equ:dBIE}) it is assumed that the boundary $\Gamma$ is sufficiently smooth.

In matrix form, this coupling reads as
\begin{equation}\label{Equ:BM}
\frac{\bphi}{2} + \bmfact \frac{\bv}{2} - \left(\bH + \bmfact \bE\right) \bphi + \left( \bG + \bmfact \bHp \right)  \bv = \bphinc + \bmfact \bv_{\text{inc}},
\end{equation}
where $\bG, \bH, \bHp,$ and $\bE$ are the matrices associated with the single layer, double layer, associated double layer, and the hypersingular operators for the Helmholtz equation. For this study, a BEM implementation  based on collocation with constant elements \cite{Kreuzer_2024} is used, so 
$$
{\bf G}_{ij} = \int\limits_{\Gamma_j} G(\bx_i,\by)\text{d}\by,
$$
the definition of the other matrices is similar. $\Gamma_j$ is a mesh element on the surface with midpoint $\bx_j$, $\bphi = (\phi(\bx_1), \dots, \phi(\bx_N))^\top$ and $\bv = (v(\bx_1), \dots, v(\bx_N))^\top$ are the vectors containing the velocity potential and the particle velocity at $\bx_i, i = 1,\dots,N$, where $N$ is the number of elements used in the discretization of the  surface $\Gamma$ of the scatterer.  $\bphi_{\text{inc}}$ and $\bv_{\text{inc}}$ are the velocity potential and the particle velocity of an external sound source, $\bmfact$ is the so-called Burton-Miller (BM) coupling factor. It has been shown \cite{BraWer65,BurMil71} that as long as the coupling factor $\bmfact$  has an imaginary part different from zero, a unique solution of the exterior problem is guaranteed for all frequencies for Neumann as well as Dirichlet boundary conditions. 

Different values for the BM factor $\bmfact$ and its sign have been discussed in the past (cf. \cite{Amini90,AmiHar90,kress85,Marburg16,Meyeretal78,terai80,Zhengetal15}), the most commonly used factor is $\bmfact = \frac{\I}{k}$, which for high frequencies is almost optimal in terms of having a small condition number for problems on the unit circle or the unit sphere. Since scattering problems involving the sphere have an analytic solution based on spherical harmonics and the spherical Bessel and Hankel functions, it is possible to derive the eigenvalues and the condition number of the linear system for the combined integral equation as a function of $\bmfact$, see also Section~\ref{Sec:Eigenvals}. Based on this analysis, it can be shown that $\bmfact = \frac{\I}{k}$ is quasi-optimal (cf.  \cite{Amini90}).

A heuristic argument for $\eta = \frac{\I}{k}$ can be given by looking at the operators involved with the dBIE Eq.~(\ref{Equ:dBIE}) in 3D: 
\begin{align*}
  H'(\bx,\by) &= \frac{\partial G}{\partial \bn_x} (\bx,\by) = G(\bx,\by) \left( \I k - \frac{1}{||\bx - \by||} \right) \frac{(\bx - \by) \cdot \bn_x}{||\bx -\by||}\\
              &= \left( O \left( \frac1r \right) + \order(k) \right) G(\bx,\by),\\
  E(\bx,\by) &= G(\bx,\by)\left[
               \left(\frac{3}{r^2} - \frac{3\I k}{r} - k^2\right) \frac{\partial r}{\partial \bn_x} \frac{\partial r}{\partial \bn_y} + \frac{1}{r}\left( \frac{1}{r} - \I k\right) \bn_x^\top \bn_y \right] = \\
              &= \left( O \left( \frac{1}{r^2} \right) + \order(k/r) + \order(k^2) \right) G(\bx,\by),
\end{align*}
where $r = || \bx - \by ||$ (cf. \cite[Eq. (15)]{Chenetal97}). %
The components of the dBIE  contain a scaling with the wavenumber $k$, at least if $k$ is big enough so that the parts involving $\order(k)$ and $\order(k^2)$ are bigger than the $\order\!\left(\frac{1}{r^2}\right)$ and $\order\!\left(\frac{1}{r}\right)$ parts. 
Using this viewpoint, multiplying the dBIE with $\frac{1}{k}$ balances the scale of BIE and dBIE. In 2D, a similar heuristic argument can be given when looking at the derivatives of the Green's function $H_0(k r)$ with respect to $r$. 

\begin{figure}[!h]
  \begin{center}
    \begin{tabular}{rcrc}
      \raisebox{0.33\textwidth}{\small a) \hspace{-20pt}} &
\includegraphics[width=0.44\textwidth]{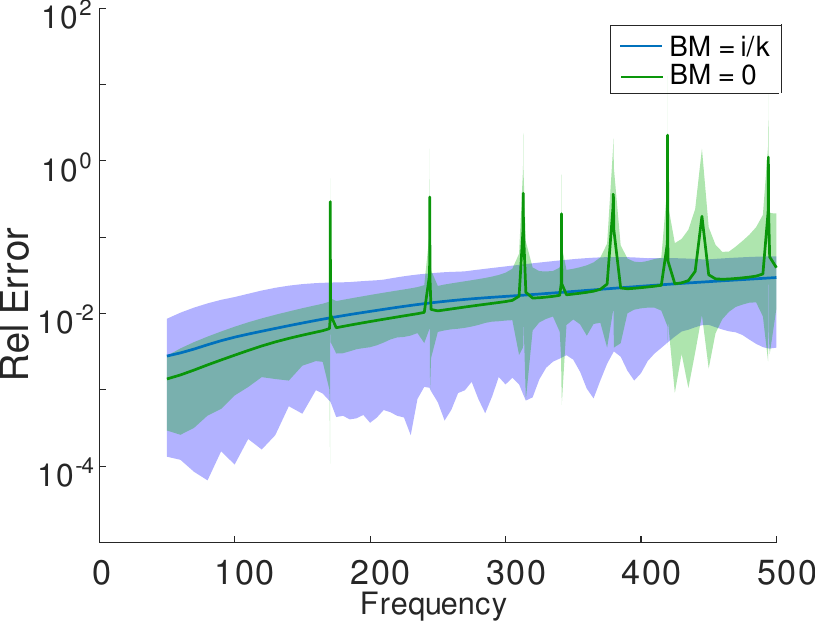}
      \raisebox{0.33\textwidth}{\small b) \hspace{-20pt}} &
    \includegraphics[width=0.44\textwidth]{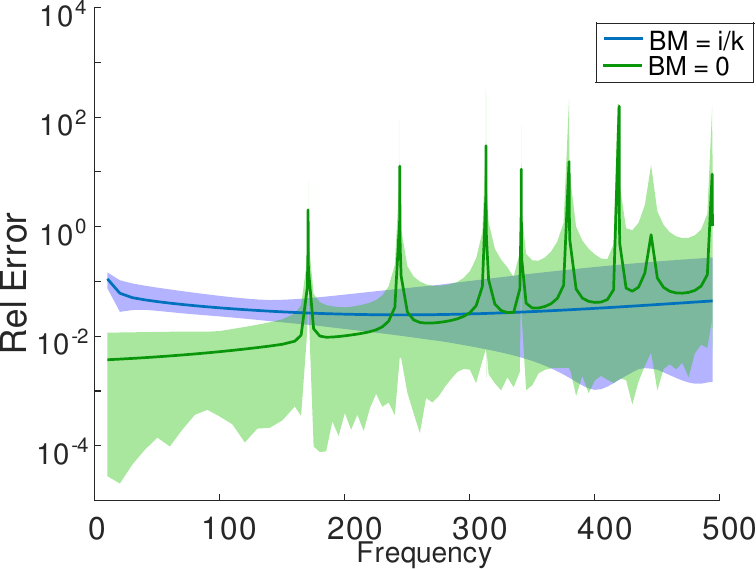}
    \end{tabular} 
  \end{center}
  \caption{Relative errors for the  problem of plane wave scattering from a) a sound hard and b) a sound soft sphere with and without Burton-Miller. The areas depict the range between maximum and minimum relative error  at the collocation nodes for  $N = 1280$ constant elements. The lines depict the mean relative error.}\label{Fig:ErrorBMnoBM}
\end{figure}

One can clearly see in Fig.~\ref{Fig:ErrorBMnoBM} that for the example of plane wave scattering from a sound hard sphere, the Burton-Miller method has a positive effect on the stability of the solution of the scattering problem. However, in the low frequency region, the relative error between numerical and analytical solution away from the critical frequencies is, in general,  smaller for the version without Burton-Miller. Using again the heuristic viewpoint, it seems that in the low frequency region, the $\order(k)$ factor becomes too small and the parts depending on the radius dominate the error. 

The data in Fig.~\ref{Fig:ErrorBMnoBM} was calculated for $N = 1280$ almost regular triangles, using collocation with constant elements \cite{Kreuzer_2024}. The highest investigated frequency of 500 Hz is slightly larger than suggested by the 6-to-8-elements per wavelength rule of thumb \cite{Marburg02}. Still, the relative error at this frequency is acceptable.

The example presented in Fig.~\ref{Fig:ErrorBMnoBM} illustrates the problem of choosing a good BM-parameter if BEM solutions at low frequencies are needed. Examples where solutions for small wavenumbers are needed are the 2.5D BEM \cite{Duhamel96,Kasess2016} or the BEM for periodic structures using the Floquet transform \cite{Clouteauetal00}. The fact that for low frequencies the coupling factor  $\eta = \frac{\I}{k}$ becomes very large, led to different approaches to define this factor for small $k$. One suggestion \cite{Amini90} is to use  $\bmfact = \min (2,\frac{1}{k}) \cdot \I$, which means for a sphere with radius $R = 1$\,m that above about 27\,Hz the BM-factor is $\eta = \frac{\I}{k}$, below it is set constant.  It was stated in \cite{AmiHar90} that ``[...] 1 has often been used in practice'', where the authors referred to the imaginary part of the coupling factor. In \cite{Duhamel96}, a factor $\bmfact = \frac{\I |k|}{\max ( |k|^2, 50)}$ was suggested in connection with the 2D and 2.5D BEM. In \cite{BruKun01}, a BM factor of $\bmfact =  \min (\frac{1}{3}, \frac{\lambda }{D})\cdot \I =   \min (\frac{1}{3}, \frac{ 2\pi}{Dk})\cdot \I$ was used\footnote{In fact, in \cite{BruKun01} the authors scale the BIE, and the factor was adapted accordingly.}, where $D$ is the diameter of the scatterer. According to \cite{BruKun01}, this factor is based on optimizing the number of necessary GMRES iterations. \cite{BruKun01} explicitly includes the diameter $D$ of the scatterer, whereas in some studies this scaling with the size of the object is implicitly assumed (cf. \cite{Amini90}: \emph{``If we denote $kR$ by $k^*$ and $\nu/R$ but $\nu^*$ then the problem is reduced to that of a unit sphere at the wave number $k^*$ with the coupling parameter $\nu^*$''}, a similar claim can be found in \cite{kress85}.).

Note that one has to be a careful when comparing different literature about the Burton-Miller method because sometimes the BIE is scaled, sometimes the dBIE. Also, when choosing the sign of the factor, it is important to look at the definitions of the harmonic time factor and the direction of the normal vector. Throughout this manuscript, it is assumed that the dBIE is scaled with the coupling factor, the harmonic time factor is set to $e^{\I \omega t}$, and the normal vectors always point away from the scatterer. Following \cite{Marburg16}, the imaginary part of the Burton-Miller factor in our numerical experiments is chosen to be positive. 

In this study, different choices of the Burton-Miller factor and their effect on the solution of a scattering problem will be investigated. As a benchmark example, the scattering of a plane wave on the unit sphere will be used. As there is some discussion about whether to adapt the BM factor to the size of the scatterer, one numerical experiment will look at the sphere with radius $R = 2$\,m. Also, a new approach based on a single step of a Modified Richardson iteration is proposed to enhance the solution of the Burton-Miller BEM in the low frequency region for sound hard problems.

In Section \ref{Sec:Eigenvals} an eigenvalue analysis of the BIE operators on the sphere is derived that helps in explaining some of the effects seen in Section \ref{Sec:NumExp}. After a short motivation for the additional iteration step in Section \ref{Sec:PostIteration}, numerical benchmark tests are done in Section \ref{Sec:NumExp}, where the solutions for different Burton Miller factors are compared with each other.  The paper is finished with a short discussion in Section \ref{Sec:Discussion}.
%%%%
%
%%%%
\section{Eigenvalue-Analysis for the sphere}\label{Sec:Eigenvals}
Following the approach proposed in \cite{Amini90}, it is possible to use the orthogonality of the Spherical Harmonics to derive the eigenvalues of the different integral operators involved with the Helmholtz equation on a sphere with radius $R = 1$\,m as function of the wavenumber $k$:
\begin{align}
  \label{Equ:EvG}\lambda_n(G) &= \I k h_n(k) j_n(k) \xrightarrow[ k \rightarrow 0]{} \frac{1}{2n+1}, \\
  \label{Equ:EvH}\lambda_n(H) = \lambda_n(H') &= \frac12 + \I k^2 h'_n(k) j_n(k) \xrightarrow[ k \rightarrow 0]{} -\frac{1}{2(2n+1)}, \\
  \label{Equ:EvE}\lambda_n(E) &= \I k^3 h'_n(k) j'_n(k) \xrightarrow[ k \rightarrow 0]{} - \frac{n(n+1)}{2n+1},
\end{align}
where $j_n, h_n, j'_n,$ and $h'_n$ are the spherical Bessel function of order $n$, the spherical Hankel function of the first kind, and their respective derivatives. From Eq.~(\ref{Equ:EvH}) it is easy to see that the hypersingular operator becomes unstable at low frequencies because one of the eigenvalues of this operator converges towards 0. A similar instability can be shown for the operator  $\frac{1}2 + H'$. Both operators are involved in the dBIE for sound hard and sound soft problems, which indicates that for small wavenumbers the dBIE for both problems is unstable. It is also easy to show that the eigenvalues for the operators $G$ and $\frac12 - H$ are bounded and sufficiently away from 0, thus, the BIE for sound hard and sound soft problems stays stable for small wavenumbers and mesh sizes $N$ that are not extensively large.\footnote{From operator theory it is known that for a linear operator $A:X \rightarrow Y$ the problem $A\phi = f$ becomes ill posed if $X$ has infinite dimension, which could negatively affect the BIE for the sound soft problem. However, for practical applications the mesh sizes are reasonable, so the ill-posedness is rather an academic problem.}

Using the asymptotic expansions \cite{AbrSte64,Amini90} of the Bessel and Neumann functions for small arguments $(k \rightarrow 0)$
\begin{align*}
  j_n(k) &\approx \frac{k^n}{1\cdot 3 \cdots (2n+1)} = \order(k^n)\\
  y_n(k) &\approx -\frac{1\cdot 3 \cdots (2n-1)}{k^{n+1}} = \order(k^{-(n+1)})
\end{align*}
it is possible to investigate the behavior of the different solutions in the low frequency region. If one looks at the eigenvalues of the operators involved in the Burton Miller approach for the Dirichlet and the Neumann problems\footnote{The values can be either found in \cite{Amini90} or directly derived by looking the values for the boundary integral operators and the fact that the Wronskian $j_n(k) h'_n(k)  - j'_n(k) h_n(k)  = \frac{i}{k^2}$.} 
\begin{align*}
  \lambda_n(D) &= \I k h_n(k) \left( j_n(k) - \eta k j'_n(k)\right),\\
  \lambda_n(N) & = \I k^2 h'_n(k) \left(j_n(k) + \eta k j'_n\right).
\end{align*}
it becomes apparent that for small $k$ and $\eta = \frac{\I}{k}$ the eigenvectors of the operators diverge:
\begin{align}
  \nonumber\lambda_n(D) & = \I k \left( \order(k^n) + \order(k^{-(n+1)})\right) \left( \order(k^n) + \order(k^{n-1})\right) = \order(k^{2n+1}) + \order(1) + \order(k^{2n}) + \order(k^{-1})\\
  \label{Equ:Operators}\lambda_n(N) & = \I k^2 \left( \order(k^{n-1}) + \order(k^{-(n+2)})\right) \left( \order(k^n) +  \order(k^{n-1})\right) = \order(k^{2n+1}) + \order(1) + \order(k^{2n}) + \order(k^{-1})
\end{align}
This explains the rise in the Burton-Miller solution in Section \ref{Sec:NumExp} towards $k = 0$ if $\bmfact = \frac{\I}{k}$. If the Burton-Miller factor scales with $k$ or is kept constant in the low frequency region, this rise in the eigenvalues can be avoided. %
%%%%%%%%%%%%%%%%%%%%%%%%%%%%%%%%%
%
%%%%%%%%%%%%%%%%%%%%%%%%%%%%%%%%%
\section{Iterative Correction of the BM Solution for Sound Hard Problems}\label{Sec:PostIteration}
It was illustrated in the previous section that in the low frequency region away from the critical frequencies, the Burton-Miller method produces larger errors than the solution of the original BIE. An alternative to using different Burton-Miller factors is presented in the following. 
As the operators and integral kernels involved with the dBIE are, in general, numerically more challenging than the operators involved with the BIE, it is plausible to assume that the error of a solution based purely on the BIE is smaller than the error for the Burton-Miller solution as long as the frequency is sufficiently far away from the set of critical frequencies, where the BIE has no unique solution (see also Fig.~\ref{Fig:ErrorBMnoBM}). This becomes also clear, when the eigenvalues for the different operators on the sphere are compared, see Section \ref{Sec:Eigenvals}. For the BIE the eigenvalues of $\frac{1}{2} - H$ converge towards $\frac{n+1}{2n+1}$ for $k \rightarrow 0$ thus they stay away from zero. For the hypersingular operator and therefore, for the dBIE,  Eq.~(\ref{Equ:EvE}) shows that for $n = 0$ the operator becomes singular. 

In \cite[Theorem 3]{Sloan76} it is shown that if a system $\phi - K \phi = 0$ with compact operator  $K$  has only the trivial solution $\phi = 0$, and if $\phi_0$ is  a Galerkin-Petrov solution of $\phi =  f + K \phi$, then $\phi_1 = f + K \phi_0$ is a better approximation of the real solution  provided that the number of ansatz functions $n$ is big enough. As the theorem is based on orthogonal projections on the ansatz and test function spaces, it cannot be applied to the collocation method directly. Also, the fact that the system needs to have a unique solution  is not fulfilled at the critical frequencies. Nevertheless, this theorem serves as a motivation for using an additional iteration step in connection with the problem of solving wave scattering around a sound hard sphere:

Let $\bphi_0$ be the solution of the Burton-Miller system Eq.~(\ref{Equ:BM}) for a sound hard scatterer using collocation with constant elements, thus
$$
\frac{\bphi_0}{2} - {\bf H}\bphi_0 = \bphi_{\inc}
$$
or
$$
\bphi_0 = 2 \bphi_{\inc} + 2 {\bf H} \bphi_0.
$$
Based on this solution, a new approximation %
\begin{equation}\label{Equ:PostIter}
  \bphi_1 = 2 \left( \bphi_{\text{inc}} + \bH \bphi_0 \right)
\end{equation}
is generated. Eq.~(\ref{Equ:PostIter}) can be interpreted as applying one single step of the method  described in \cite{Sloan76}, where the operator $K$ is replaced by the discretized version of the operator $2H$ used in the BIE. Alternatively, the approach can be interpreted as one step of a modified Richardson Iteration/Steepest Descent method 
$$
\bphi_{n+1} = \bphi_n + w \left( \bphi_{\text{inc}} - \left(\frac{\bf I}{2} - {\bf H}\right) \bphi_n\right)
$$
with step size $w = 2$ applied to the discretized BIE, where ${\bf I}$ denotes the identity matrix of size $N\times N$. Note that, in general, the step size of $w = 2$ is too large to guarantee convergence of the modified Richardson Iteration\footnote{The condition for convergence is $||{\bf I} - w ({\bf I}/2 - {\bf H})||<1$. }. For the example of plane wave scattering on a sound hard unit sphere at $f = 50$\,Hz, the step size should be $w < 1.6$ to guarantee convergence, at $f = 200$\,Hz, the step size should theoretically be $w < 1.3$. However, the numerical experiments in Section \ref{Sec:NumExp} show, that one single step with $w = 2$ can  ``correct'' the Burton-Miller solution towards the solution of the BIE \emph{without} the singularities around the critical frequencies. It is possible to use a smaller step size, but in that case, one single step may not be enough to reach ``optimal'' accuracy.

Note that the approach described in \cite{Sloan76} is fitted to equations of the type $\phi = f + A\phi$, where $A$ is some linear operator. This is also a restriction of the approach proposed here that prohibits an efficient solution of a sound soft problem by this approach. Theoretically, it is possible to use the dBIE 
$$
\frac{\bv}{2} + {\bf H}' \bv = \bv_\inc
$$
to construct a similar iteration
$$
{\bf v}_{n+1} = 2 \left( {\bf v}_{\text{inc}} - {\bf H}' {\bf v}_n \right),
$$
however, the discussion in Section \ref{Sec:Eigenvals} with respect to the eigenvalues of the different boundary integral operators shows that for the example of the sphere one of the eigenvalues for the associated double layer potential $H'$ is $-1/2$, see Eq.~(\ref{Equ:EvH}). Thus the stability of the dBIE deteriorates for small $k$ and an additional iteration step just based on this equation will not lower the error. 

\section{Numerical Experiments}\label{Sec:NumExp}
\subsection{Sound Hard  Sphere}
As the first benchmark example, the standard problem of plane wave scattered at a sound hard sphere is used. This example has the advantage, that an analytic
solution of the problem exists \cite[Eq.~(6.186)]{Williams99}, and it is therefore possible
to easily calculate the numerical error of the BEM method. In all numerical experiments, the sphere was discretized using a mesh with almost regular triangles based on a discretization of a projected subdivided icosaeder. Two meshes with $N = 1280$ and $N = 5120$ elements were used (see also Fig.~\ref{Fig:Meshes}), the radii of the sphere were set to $R = 1$\,m and $R = 2\,$\,m, and the linear systems for all benchmark problems will be solved using an LU-decomposition implemented in LAPACK \cite{Andersonetal99}. In all computations it is assumed that the speed of sound and the density are given by $c = 340$\,m\,s$^{-1}$ and $\rho = 1.3$ kg\,m$^{-3}$.
\subsubsection{Sphere with radius $R = 1$\,m}
\begin{figure}[!h]
  \begin{center}
    \begin{tabular}{rcrc}
      \raisebox{0.4\textwidth}{\small a) \hspace{-16pt}} &                                                            \includegraphics[width=0.44\textwidth]{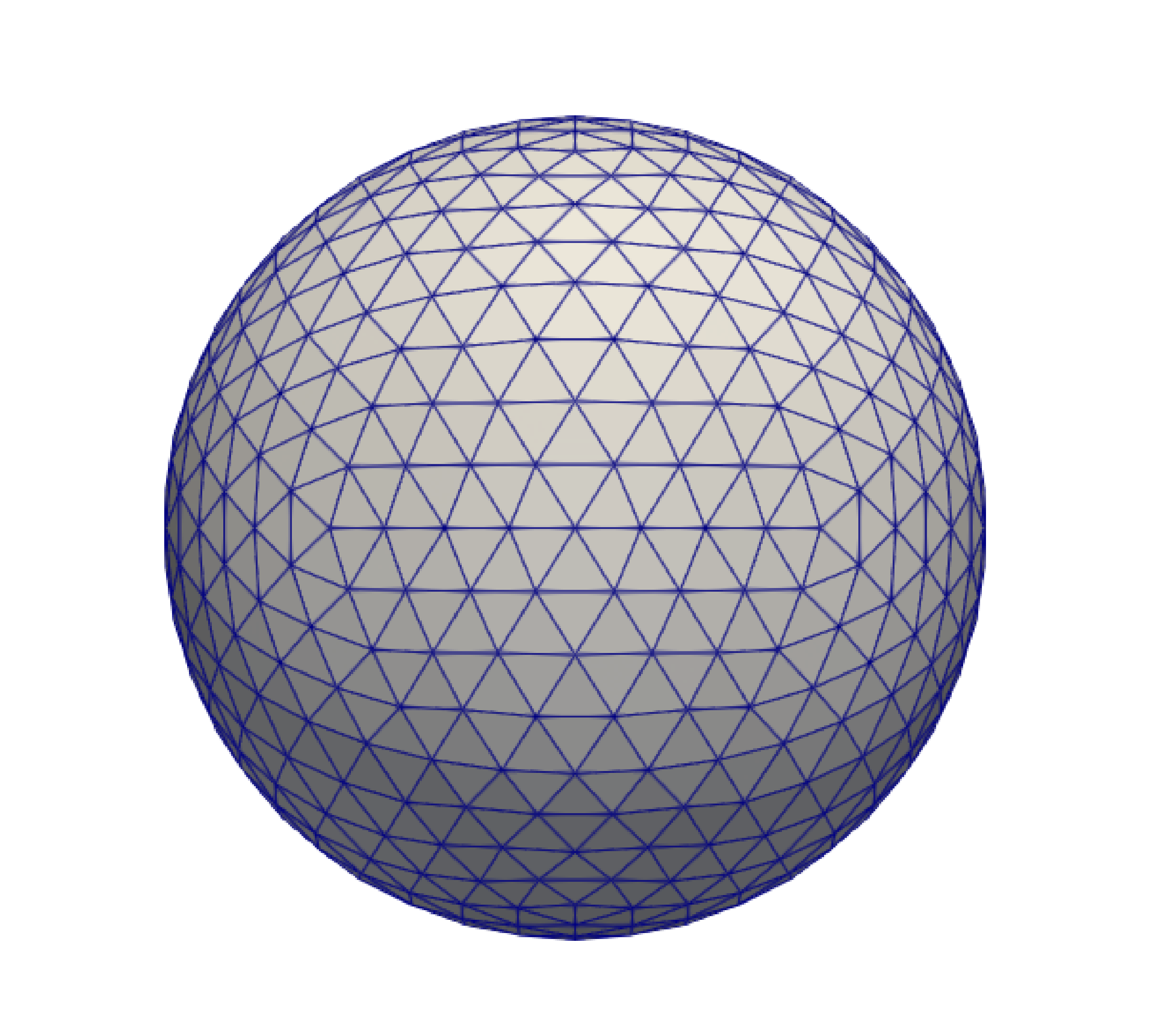}
      &
        \raisebox{0.4\textwidth}{\small b) \hspace{-16pt}}
      &
        \includegraphics[width=0.44\textwidth]{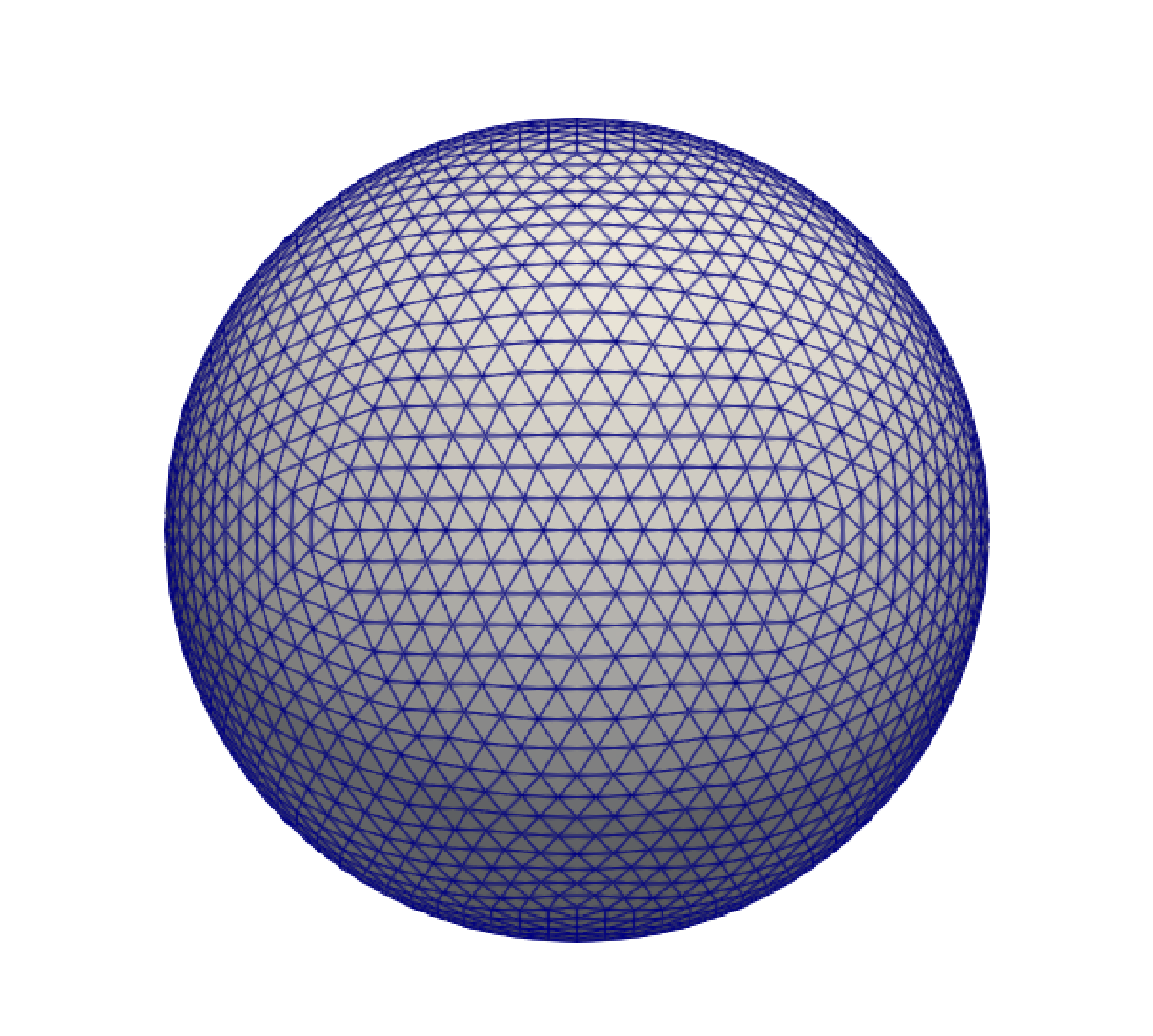}
    \end{tabular}
    \caption{Discretization of the unit sphere with a) $N = 1280$ and b) $N = 5120$ almost regular triangles.}\label{Fig:Meshes}
  \end{center}
\end{figure}

In the first experiments, the condition numbers and the relative errors of the BEM for different BM factors will be investigated for the example of plane wave scattering from a sphere with radius $R = 1$\,m that was discretized with $N = 1280$ and $N = 5120$  triangular elements (see also Fig.~\ref{Fig:Meshes}). The frequencies of interest were between $10$\,Hz and $500$\,Hz in steps of 10\,Hz. For a frequency of 500 Hz, the average edge length of the elements of the mesh was slightly over the limit recommended by the 6-to-8 elements per wavelength rule of thumb \cite{Marburg02}.  Around the critical frequencies, which are given by the zeros of the spherical Bessel function $j_n(k)$ and which are summarized in Tab.~\ref{Tab:Crit1280}, 9 steps with step size $0.05$\,Hz were made. Because of the discretization of the sphere's surface and the fact that the solution is calculated at the midpoints of the elements, the critical frequencies for the BEM calculation are not exactly at the zeros of the Bessel functions. These frequencies were determined by finding the frequency with a large condition number close to the theoretical critical frequency. These values are close to the theoretical frequencies for spheres where the radii are determined by the mean norm of the collocation nodes.

\begin{figure}[!h]
  \begin{center}
    \begin{tabular}{lclc}
    \raisebox{0.35\textwidth}{{\small a)}\hspace{-10pt}} &\includegraphics[width=0.44\textwidth]{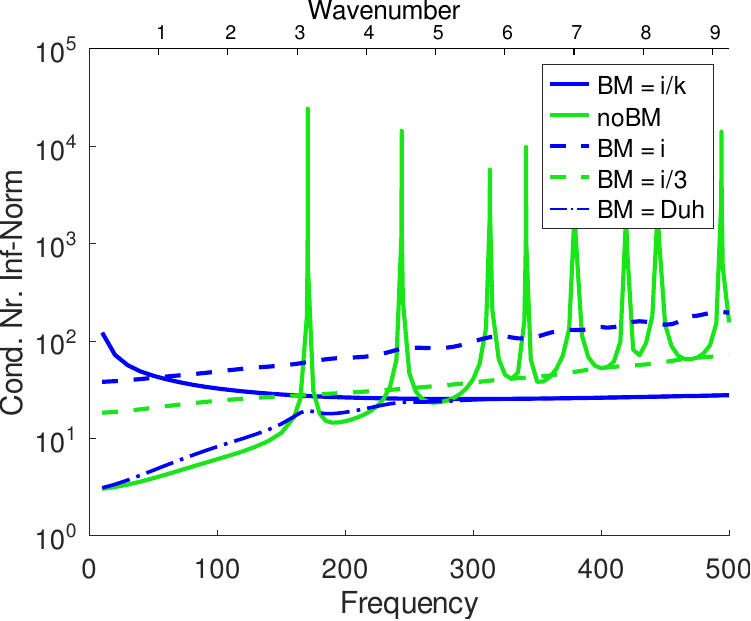} &
    \raisebox{0.35\textwidth}{{\small b)}\hspace{-10pt}} &\includegraphics[width=0.43\textwidth]{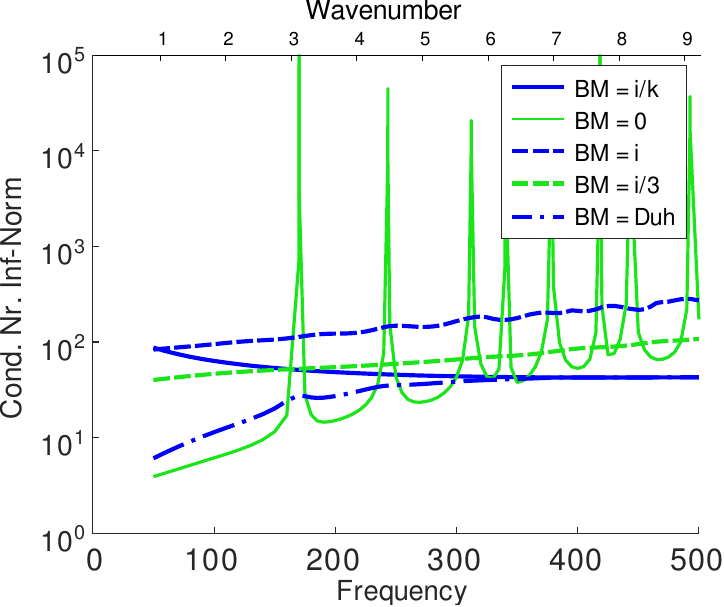}
    \end{tabular}
     \caption{Condition numbers for different Burton-Miller coupling factors  $\bmfact = 0,\frac{\I}{k}, \frac{\I}{3}, \I$ and $\bmfact$ chosen according to \cite{Duhamel96}. For the discretization of the unit sphere a) N = 1280, b) N = 5120 triangular elements were used.}\label{Fig:CondandRes}
  \end{center}
\end{figure}

In Fig.~\ref{Fig:CondandRes}, the condition numbers of the BEM system matrices ${\bf I}/2 - {\bf H} - \bmfact {\bf E}$ (see also Eq.~(\ref{Equ:BM})) for different Burton-Miller factors ($\bmfact = \frac{\I}{k}, 0, \I, \frac{\I}{3}$ and the factor proposed by Duhamel in  \cite{Duhamel96}) are depicted as function of the frequency for $N = 1280$. Note that for the frequency range of interest, $\min (1/3, \frac{2\pi}{2rk}) = 1/3$, thus the factor proposed in \cite{BruKun01} is reduced to $\bmfact = \frac{\I}{3}$.
\begin{table}
  \caption{Rounded theoretical and numerically determined (up to three digits) critical frequencies for the sound hard unit sphere  with $N = 1280$ and $N = 5120$ elements. The numerically critical frequencies were determined by looking at the condition number of the BIE.}\label{Tab:Crit1280}
  \begin{center}
    \begin{tabular}{ccccc}
      \hline
      Theor. & 170.00 & 243.150 & 311.876 & 340.000 \\ 
      $N = 1280$  & 170.656 & 244.083 & 313.062 & 341.315\\ 
      $N = 5120$  & 170.164 & 243.383 &  312.172 & 340.3280 \\
      \hline
    \end{tabular}
    \begin{tabular}{ccccc}
      Theor. & 378.136 & 418.034 & 442.781 & 492.155\\
      $N = 1280$ & 379.586 & 419.65 & 444.445 & 494.045\\
      $N = 5120$ & 378.498 & 418.436 &443.195 & 492.627\\
    \end{tabular}
  \end{center}
\end{table}

\begin{figure}[!h]
  \begin{center}
    \includegraphics[width=0.60\textwidth]{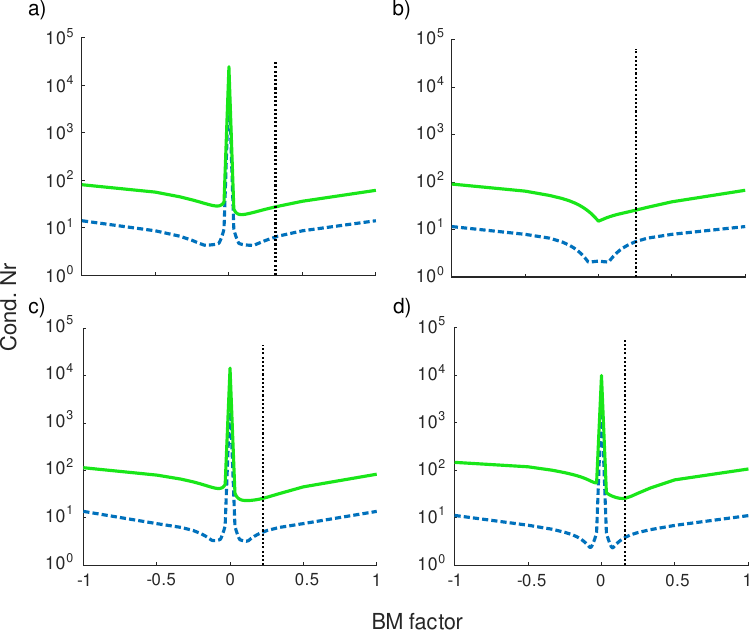}
  \end{center}
  \caption{Condition numbers of the linear system for N = 1280 at a) 170.656 Hz, b) 200 Hz, c) 244.083 Hz and d) 343.12 Hz as a function of the Burton-Miller factor in the interval $[-1,1]$. The (blue) dashed curves depict the condition number based on the 2-norm, the (green) continuous curves the condition number based on the Inf-norm. The dotted lines represent $\eta = \frac{\I}{k}$.}\label{Fig:CondNrs}
\end{figure}

The condition numbers were calculated using the Inf-Norm  by the LAPACK-Routine \texttt{zgecon} \cite{Andersonetal99}. This condition number is about $N$ times larger than the ``usual'' condition number with respect to the 2-norm, see, e.g., Fig.~\ref{Fig:CondNrs} for a comparison of the two different definitions. Up to 250\,Hz the condition number for $\bmfact = 0$ is smaller than the condition number for the versions with $\bmfact \neq 0$, except close to the critical frequencies. In general, the choice of the factor according to \cite{Duhamel96} (represented by the dash-dotted line) is very close to the ``optimal'' values, also the switching between $\bmfact = \frac{\I k}{50}$ and $\bmfact = \frac{\I}{k}$ is quite accurate to keep the overall condition number small. It is also clearly visible that the factor $1/k$ has a negative effect on the condition number for very low frequencies. However, the all values for the condition numbers are in a feasible range. 

In Fig.~\ref{Fig:CondNrs}, the condition numbers at $f = 170.656$, $200, 244.083$, and $343.315$\,Hz are depicted as function of the Burton-Miller factor $\eta$ for $N = 1280$ elements. Especially for low frequencies, $\eta = \frac{\I}{k}$ (denoted by the dotted line) does not lead to the smallest condition number. This may be because a) of the geometric error introduced by the mesh, which means that the solution of the BEM is different from the analytic solution, b) the frequencies investigated were not the exact critical frequencies and c)  that the wavenumber dependent parts of the derivatives of the Green's function were too small compared to the $\order(\frac{1}{r})$ and $\order(\frac{1}{r^2})$ factors. In general, most papers dealing with finding a quasi optimal Burton-Miller factor look at the unit sphere and the \emph{analytic} solution and eigenvalues for the non-discretized sphere, respectively. As relatively few elements are used for the BEM, there is already some difference in geometry to the ``optimal'' case. Fig.~\ref{Fig:CondNrs}b, where a non-critical frequency was investigated, also shows that  the non Burton-Miller version $\bmfact = 0$ is better conditioned for the sound hard sphere. 
\begin{figure}[!h]
  \begin{center}
     \includegraphics[width=0.46\textwidth]{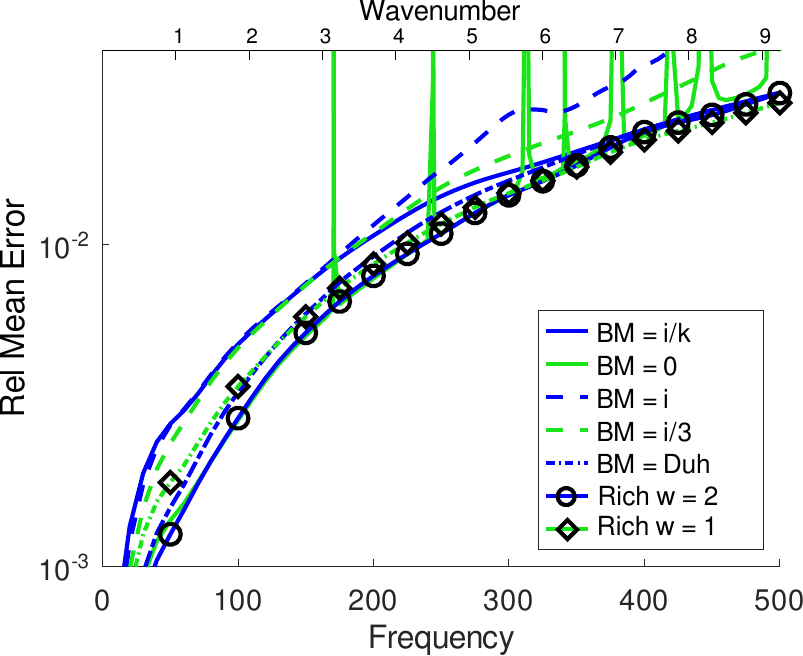}
     \includegraphics[width=0.46\textwidth]{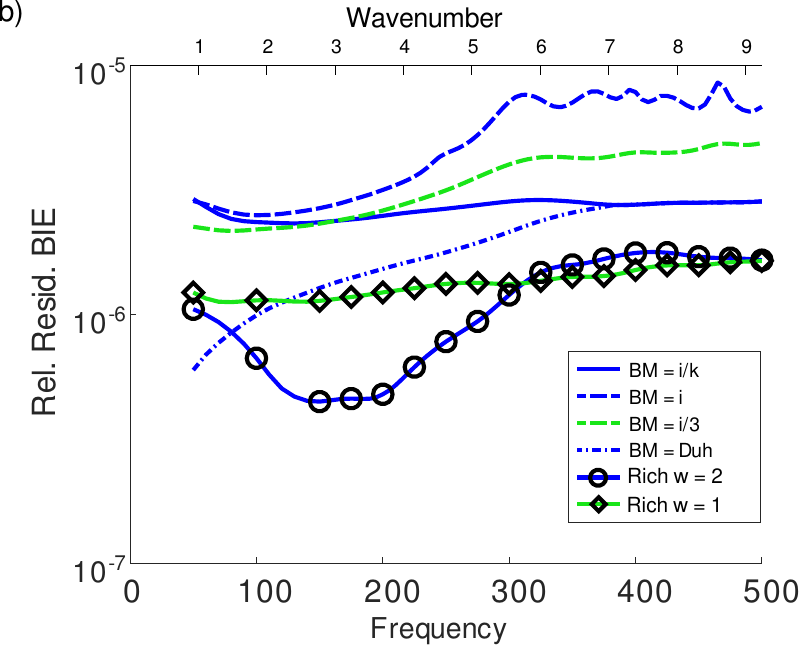}
     \caption{a) Relative errors of the  BEM solutions compared to the  analytic solution on the sound hard sphere for the different Burton-Miller factors (lines) and the relative errors after a single modified Richardson iteration step (symbols) with $w = 2$ and $w = 1$. 
       b) Residuum generated by inserting the Burton-Miller solutions into the BIE. The number of elements used in the mesh is $N = 1280$. }\label{Fig:MeanRelErrors1280}
  \end{center}
\end{figure}

\begin{figure}[!h]
  \begin{center}
    \includegraphics[width=0.46\textwidth]{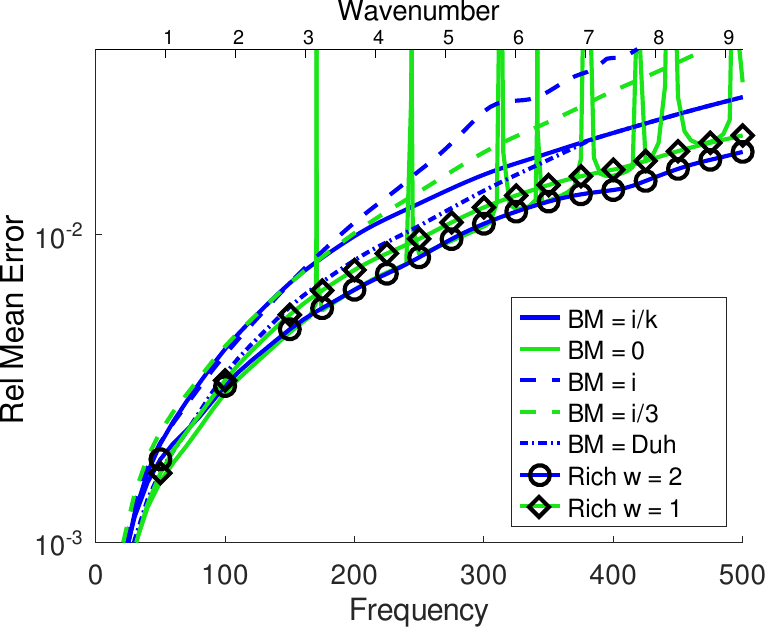}
    \caption{Mean relative errors for the sound hard sphere where the analytic solution is calculated at a sphere with $r \approx 0.996$\,m.} \label{Fig:MeanRelErrors1280NoUnit}
  \end{center}
\end{figure}

Naturally, the different condition numbers have an impact on the relative error between BEM solution and analytical solution.  
Fig.~\ref{Fig:MeanRelErrors1280}a depicts the mean relative error over all collocation nodes where the analytic solution is determined on the unit sphere, thus the comparison  also includes the geometric error caused by discretization of the sphere and the fact that the collocation nodes do not lie on the unit sphere.  Again, the error is almost optimal if $\bmfact$ is chosen according to \cite{Duhamel96} (dash-dotted line). Additionally, the mean relative errors of one single step of the modified Richardson iteration (see Section~\ref{Sec:PostIteration}) with $w = 2$ (circles) and $w = 1$ (diamond symbols) are depicted. As one can see, one single step with $w = 2$ reduced the relative error for the Burton-Miller formulation to match the error level of the original BIE formulation, at least in the low frequency region. Above 400 Hz, the error for $w = 2$ gets slightly worse than the error for $\eta = \frac{\I}{k}$ and the error for the Richardson step with $w = 1$.

It was already mentioned that for this error plot the analytic solution is calculated for points on the unit sphere. For $N = 1280$ elements, the average radius over all collocation node is about $r_c \approx 0.996$\,m. If the analytic solution is calculated at the collocation nodes (or for a sound hard sphere with radius $R = 0.996$\,m) the behavior of the mean error is slightly different, see Fig.~\ref{Fig:MeanRelErrors1280NoUnit}. In this setting, the difference between the BEM solution with Burton-Miller ($\bmfact = \frac{\I}{k}$) and without Burton-Miller becomes smaller in the low frequency region, but it gets bigger for larger frequencies. Here, the additional modified Richardson step results in the smallest error.

\begin{figure}[!h]
  \begin{center}
    \begin{tabular}{lclc}
      \raisebox{0.35\textwidth}{ {\small a)}\hspace{-10pt} } & \includegraphics[width=0.42\textwidth]{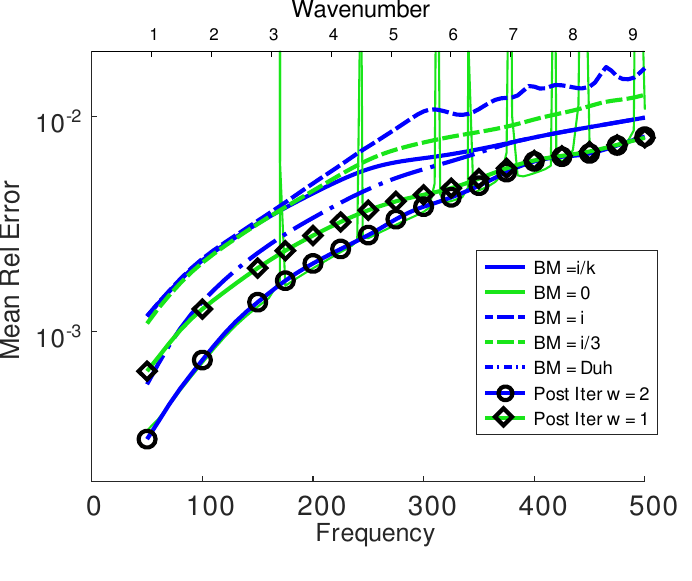} &
                                                                             \raisebox{0.35\textwidth}{\small b)\hspace{-10pt}} & \includegraphics[width=0.42\textwidth]{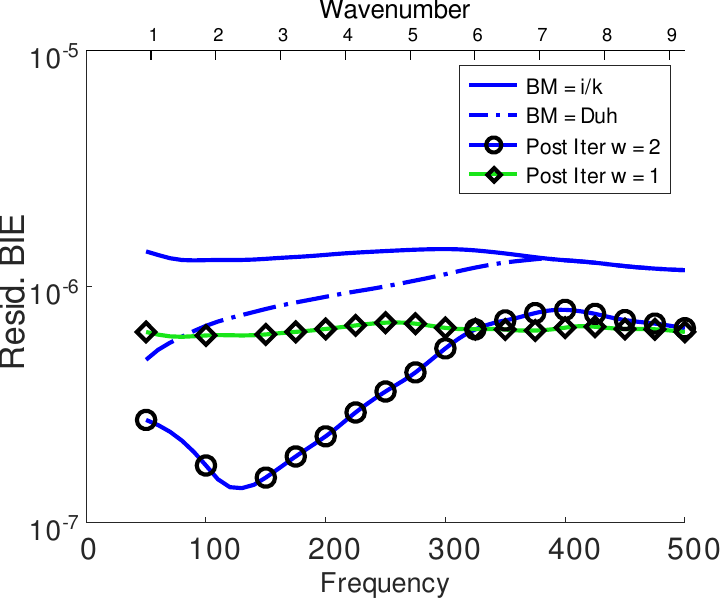}
    \end{tabular}
    \caption{a) Relative errors of the BEM solutions compared to the  analytic solution at the collocation nodes on the sphere for the different Burton-Miller factors. The lines depict the mean relative error  over all collocation nodes. b) Residuum generated by inserting the Burton-Miller solution into the BIE. The number of elements in the mesh is $N = 5120$. }\label{Fig:MeanRelErrors5120}
  \end{center}
\end{figure}

On a sidenote, the fact that a solution $\bphi_{\text{BM}}$ solves the discretized combined boundary integral equation with an relative residuum of about $10^{-16}$, does not necessarily mean, that this solution also solves the original discretized BIE with a similar accuracy. 
Fig.~\ref{Fig:MeanRelErrors1280}b illustrates this behavior of the Burton Miller solution that was calculated using a direct solver from LAPACK. For $f = 300$~Hz, for example, the mean residuum of the Burton Miller solution $\bphi_{\text{BM}}$ inserted in the discretized BIE 
$$
\beps_{\text{BIE}} = \left(\frac{\bf I}{2} - {\bf H} \right) \bphi_{\text{BM}} - \bphi_\inc
$$
is on average about $3 \cdot 10^{-6}$, a similar range can be observed if $\bphi_{\text{BM}}$ is inserted in the dBIE:
$$
\beps_{\text{dBIE}} = -{\bf E} \bphi_{\text{BM}} - \bv_\inc.
$$
Only the  combination $\beps_{\text{BIE}} + \eta \beps_{\text{dBIE}} \approx 10^{-16}$ is in the range of the machine precision. Although not really relevant for the real error of the solution compared to the analytic solution, it is, at least in the opinion of the author, an interesting behavior of the Burton Miller method. 

In Fig.~\ref{Fig:MeanRelErrors5120}, similar data is plotted for a mesh with $N = 5120$ elements in the same frequency range. In the numerical experiment, the distance between the Burton-Miller solution and the solution of the BIE gets bigger. In this case, the factor chosen according to Duhamel still yields an increase in accuracy compared to $\eta = \frac{\I}{k}$, but the additional Richardson iteration step with $w = 2$ yields an even smaller mean error, while still being stable at the critical frequencies. In Fig.~\ref{Fig:MeanRelErrors5120}b, the residuum of the BM-solution with $\eta = \frac{\I}{k}$, $\eta$ chosen according to Duhamel, and the solution of the  iteration step with $w = 2$ and $w = 1$ inserted into the BIE is depicted.

\subsection{Sphere with $R = 2$\,m, $N = 5120$}
As the BM-factor is often derived by investigating the unit sphere as the scatterer (cf. \cite{Amini90,AmiHar90,kress85,terai80,Zhengetal15}), the influence of the size of the scatterer will be investigated in this subsection. In the following, condition numbers, residuum, and relative errors for the ``classical'' BM-factor $\bmfact = \frac{\I}{k}$, for $\bmfact = \frac{\I}{2k}$, the factors proposed by \cite{BruKun01} and \cite{Duhamel96}, and the data for the iteration approach will be investigated for the example of a sound hard sphere with radius $R = 2$\,m.
\begin{figure}
  \includegraphics[height=0.4\textwidth]{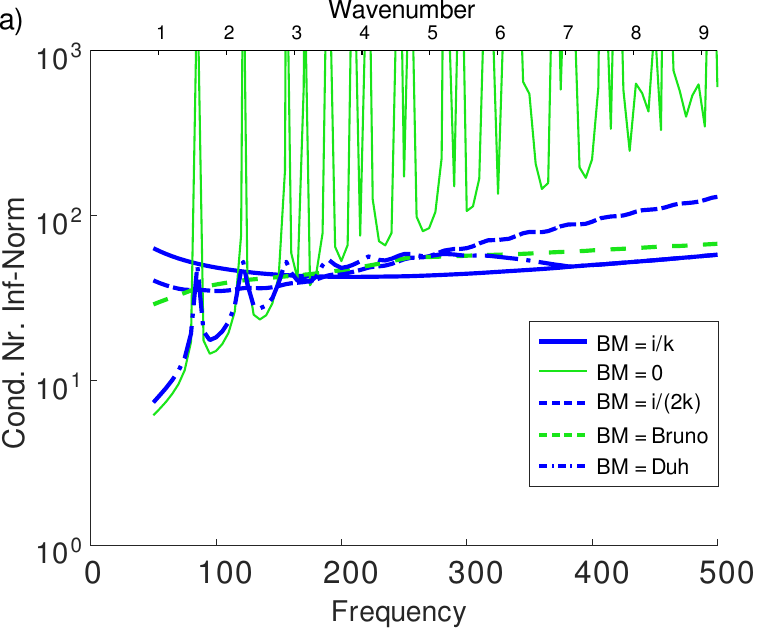}
  \includegraphics[height=0.4\textwidth]{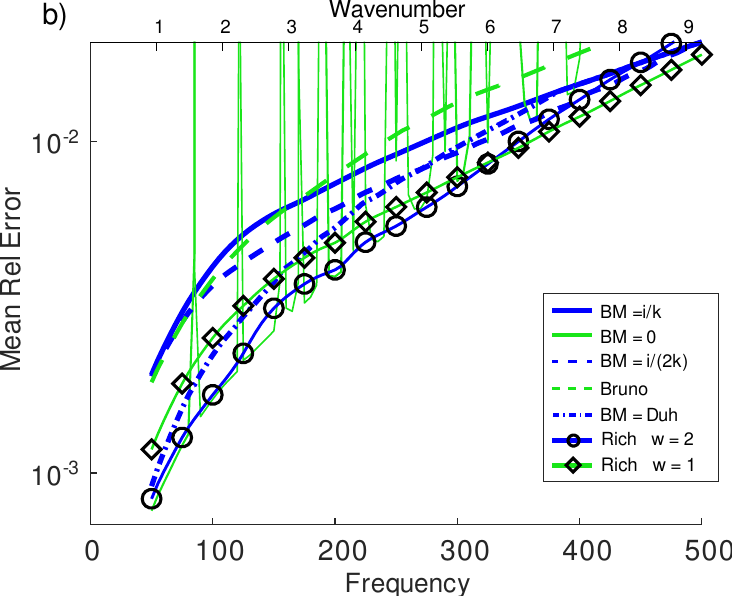}
  \caption{a) Condition numbers and b) mean relative error for different choices of the BM factor for the sphere with radius $R = 2$\,m. The sphere was discretized using $N = 5120$ elements.}\label{Fig:ErrCondr2}
\end{figure}
It is clear that for $R = 2$\,m the number of critical frequencies in the observed frequency range grows as they now depend on the roots of $j_\ell(2k),\ell = 0,1,\dots$.  Thus, the condition number of the BEM without Burton-Miller grows earlier illustrating the need for having the regularization provided by the BM-method.

For $R = 2$\,m  and $k > \pi$ the condition number for the non-Burton-Miller version is always higher than the number for the versions with Burton-Miller for all investigated choices of $\bmfact$ (see for example Fig.~\ref{Fig:ErrCondr2}). The influence of the critical frequencies on the BIE becomes dominant for  $k > \pi$, and if the BM-factor is interpreted as a weighting of BIE and dBIE, it is reasonable that the condition number for  $\eta = \frac{\I}{2k}$  is bigger than for $\eta = \frac{\I}{k}$. %

When comparing the results for $\bmfact = \frac{\I}{k}$ and $\bmfact = \frac{\I}{2k}$ it becomes clear that for relatively low frequencies including a BM factor that is dependent on the size of the scatterer has an advantage. However, for higher frequencies the non-scaled factor $\bmfact = \frac{i}{k}$ yields better mean errors. For frequencies $f = 600$\,Hz  and $f = 800$\,Hz the mean relative errors for $\bmfact = \frac{\I}{k}$ are about $\epsilon_{\text{mean}} = 0.026$ and $\epsilon_{\text{mean}} = 0.049$, for $\bmfact = \frac{\I}{2k}$ the errors are about $\epsilon_{\text{mean}} = 0.029$ and $\epsilon_{\text{mean}} = 0.061$.

Again, the factor proposed by \cite{Duhamel96} seems to be a very good choice for the BM factor. 
If the BM solution is modified using a single Richardson iteration step the error decreases in the low frequency region for $w = 2$ as well as $w = 1$. For the higher frequencies  the step size $w = 2$ is too big, in this case the iteration with step size $w = 1$ yields better results.

The analytic solution was calculated based on the collocation nodes, however  the discretization with $N = 5120$ is fine enough and there is only little difference between the analytic solution on the sphere and the analytic solution at the collocations nodes, which are slightly inside the sphere.  The version using the BM factor suggest by \cite{BruKun01} has a slightly smaller error in this case.
\subsection{Sound Soft Sphere}
The approach proposed by Burton and Miller can be ``readily extended to allow for complex values of $k$ and more general boundary conditions''~\cite{BurMil71}. It was also shown in \cite{BraWer65,Kussmaul69} that the non-uniqueness problem in connection with Dirichlet boundary conditions can be avoided by using a suitable combination of single- and double-layer potentials.

\begin{figure}[!h]
  \begin{center}
    \begin{tabular}{lclc}
      \raisebox{0.35\textwidth}{{\small a)}\hspace{-10pt}} &
      \includegraphics[width=0.42\textwidth]{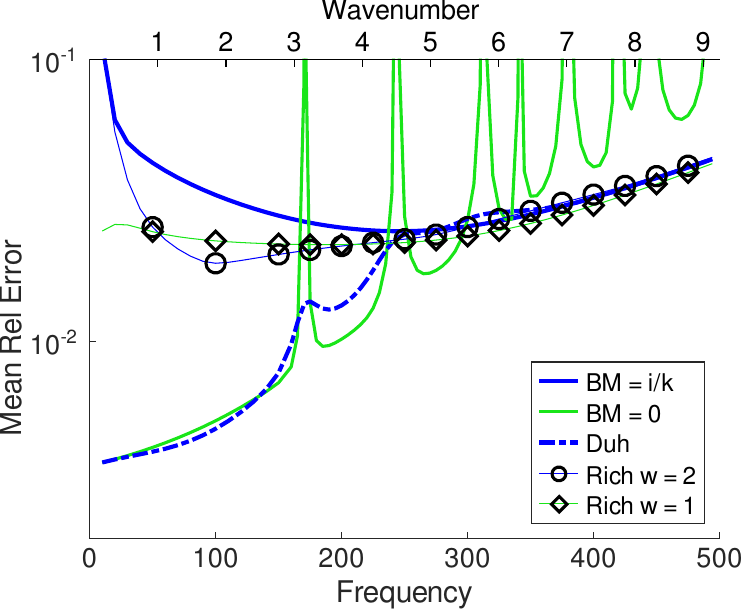} &
      \raisebox{0.35\textwidth}{{\small b)}\hspace{-10pt}} &
      \includegraphics[width=0.42\textwidth]{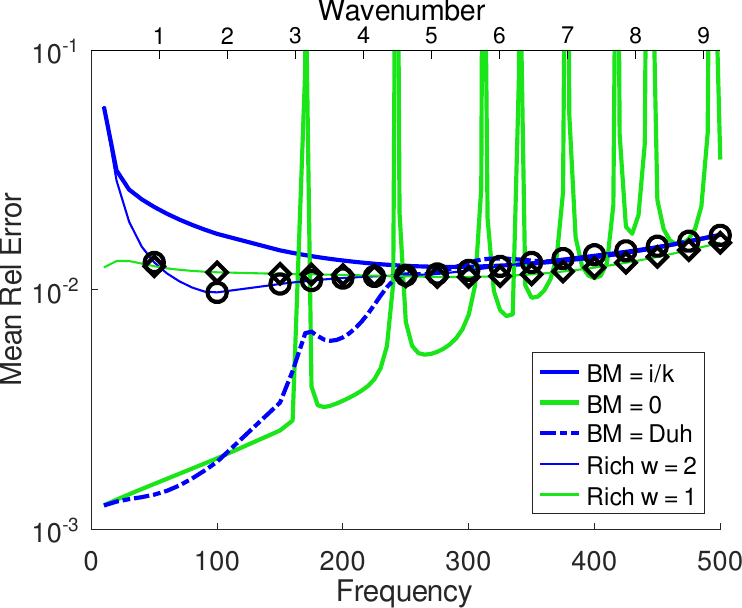}
    \end{tabular}
    \caption{Relative errors of the  BEM solutions compared to the  analytic solution on the sound soft sphere for the different Burton-Miller factors (lines) and the relative errors after a single Richardson iteration step (symbols) with $w = 2$ and $w = 1$ for a) N = 1280 elements and b) N = 5120 elements.}\label{Fig:ErrorSoft}
  \end{center}
\end{figure}

For a sound soft sphere with $\phi(\bx) = 0$ for  $\bx \in \Gamma$, the difference between BM-solution and the solution of the BIE becomes larger for smaller frequencies (cf. Fig.~\ref{Fig:ErrorSoft}). In fact, the solutions of the dBIE
\begin{equation}\label{Equ:SoftdBIE}
  \frac{\bv}{2} + {\bf H'}\bv = \bv_{\inc}
\end{equation}
becomes unstable for $k \rightarrow 0$ because one of the eigenvalues of $\bf H'$ converges towards $-\frac{1}{2}$ for $k \rightarrow 0$. This explains the sharp rise in the error of the BM solution towards $f = 0$.
Also, it is not trivial to construct one single and ``cheap'' iteration step to increase the accuracy of the BM solution anymore.  As already mentioned in Section \ref{Sec:PostIteration} an iteration based on the dBIE is not recommended because of stability problems. However, it is possible to use one or two steps of a Richardson Iteration applied to the BIE
$$
{\bf G}\bv = \bphi_{\inc},
$$
but for the reason mentioned above  the Burton Miller solution is not a good starting point for small $k$ and one single Richardson Iteration step does not yield a substantial gain in accuracy as Fig.~\ref{Fig:ErrorSoft} illustrates.

As in all other cases, the approach proposed by \cite{Duhamel96} seems to be a good compromise between accurate and regular solutions, although it  does partly not yield the ``optimal'' BM-factor.
%%%%%%%%%%%%%%%%%%%%%%%%%%%%%%
% 
%%%%%%%%%%%%%%%%%%%%%%%%%%%% 
\section{Discussion}\label{Sec:Discussion}
Regularization methods like the Burton-Miller method play an important role for the BEM in acoustics. In general, the frequency dependent coupling factor $\eta = \frac{\I}{k}$ is a good choice for coupling BIE and dBIE to regularize the solution of the BIE for most frequencies. However, the numerical experiments made in the previous sections showed that  especially in the low frequency region, where the critical frequencies are sparse and the numerical instabilities of the BIE are concentrated around  those frequencies, a different choice for the BM factor can have a positive effect on the condition number and in turn on the calculation error. Originally, the coupling factor $\eta = \frac{\I}{k}$ has been derived by investigating the condition number of scattering problems involving the unit sphere for relatively high frequencies, but there is always the difference between numerical approximation using, e.g., the BEM, and the analytic representation of the solution using spherical harmonics. This difference can be seen, for example, when comparing the critical frequencies of the analytic approach with the frequencies with high condition numbers for the numerical approximation.

When looking at Fig.~\ref{Fig:CondandRes}, it becomes clear that apart from the critical frequencies, the condition numbers in the low frequency region for the BIE is smaller than for the Burton-Miller formulation. From the eigenvalue analysis in Section \ref{Sec:Eigenvals} it becomes clear that for the problem on the sphere the dBIE and thus the Burton-Miller system becomes unstable if $k \rightarrow 0$, which explains the smaller condition number for the version without using the Burton-Miller method.

Also, the frequency range where the condition number is, in general, smaller for the BIE than for the BM formulation grows with the number of elements. For a sound hard sphere, the Burton-Miller method involves the hypersingular operator, which is numerically challenging. It is reasonable to assume that the numerical errors introduced by the discretization and the calculation of this operator have a negative effect on the condition number and the accuracy if the coupling factor becomes large in the low frequency region. For a sound soft sphere, the accuracy of the numerical solution of the dBIE gets small for $f \rightarrow 0$ because the system itself becomes singular. In this case, a solution purely based on the BIE is essential.

Using a BM factor in the low frequency region that is \emph{not} dependent on $\frac{1}{k}$   as, e.g., proposed by \cite{Duhamel96} or for larger objects the factor dependent on the cross-section of the scatterer as proposed by \cite{BruKun01}  has a positive effect on the condition number and the accuracy of the solution as also illustrated by the eigenvalue analysis for a sphere in Section \ref{Sec:Eigenvals}. The factors proposed by \cite{BruKun01, Duhamel96} have the slight drawback, that their definition changes at a certain frequency, and it is not clear, how this frequency is motivated. In \cite{Duhamel96} this frequency was defined by $\max( |k|^2, 50)$, which works excellent for the  sphere of radius $R = 1$\,m (cf. Fig.~\ref{Fig:CondandRes}), but leads to higher than necessary condition numbers for a bigger sphere (see for example Fig.~\ref{Fig:ErrCondr2}). Nevertheless, as a one-glove-fits-it-all BM factor the factor proposed in \cite{Duhamel96} is a very good choice. Compared to the usual choice $\bmfact = \frac{\I}{k}$ it has the advantage that the systems for both sound hard and sound soft problems stay stable for small $k$ (see Eq.~(\ref{Equ:Operators})). 

As an alternative to using different BM factors for different frequency regions, it was also investigated how one additional step of a modified Richardson iteration for the BIE using the BM-solution with $\eta = \frac{\I}{k}$ as a starting value increases the accuracy of the BM solution. This approach is computationally relatively cheap, as it involves only one additional matrix-vector multiplication plus two vector additions. Of course, one may argue, that the additional iteration step may be unnecessary for higher frequencies, and that it would be hard to define beforehand when this additional iteration should be made or not. On the other hand, it would be easy to compare the solution vector before the additional iteration step and after. If the difference is below a certain tolerance, one could cease to use the additional iteration for frequencies higher than the current one. However, this additional iteration step only yields low errors for sound hard boundary conditions. 

\section{Acknowledgements}
The author would like to thank the reviewers for their constructive remarks, which helped to improve the paper a lot.

\bibliography{BMpaper}
\bibliographystyle{plain} 
\end{document}